\documentclass[11pt]{article}
\usepackage{amssymb,latexsym}
\usepackage{epsfig}
\usepackage{eufrak}
\usepackage{amsmath}
\usepackage{mathrsfs}
\usepackage{color}

\newtheorem{theorem}{Theorem}
\newtheorem{lemma}{Lemma}
\newtheorem{corollary}{Corollary}
\newtheorem{remark}{Remark}

\newcommand{\be}{\begin{equation}}
\newcommand{\ee}{\end{equation}}
\newcommand{\bee}{\begin{eqnarray*}}
\newcommand{\eee}{\end{eqnarray*}}
\newcommand{\bel}{\begin{eqnarray}}
\newcommand{\eel}{\end{eqnarray}}
\newcommand{\bec}{\begin{cases}}
\newcommand{\eec}{\end{cases}}
\newcommand{\bem}{\begin{bmatrix}}
\newcommand{\eem}{\end{bmatrix}}

\newcommand{\bed}{\begin{description}}
\newcommand{\eed}{\end{description}}
\newcommand{\bei}{\begin{itemize}}
\newcommand{\eei}{\end{itemize}}
\newcommand{\ben}{\begin{enumerate}}
\newcommand{\een}{\end{enumerate}}

\newcommand{\beL}{\begin{lemma}}
\newcommand{\eeL}{\end{lemma}}
\newcommand{\beT}{\begin{theorem}}
\newcommand{\eeT}{\end{theorem}}

\newcommand{\bpf}{\begin{pf}}
\newcommand{\epf}{\end{pf}}

\newcommand{\pfbox}{\hfill\mbox{$\Box$}}

\newenvironment{pf}{\paragraph*{Proof{\rm.}}}{\pfbox\bigskip}

\begin{document}

\title{{\bf Constrained Optimal Synthesis and Robustness Analysis
 by  Randomized Algorithms\thanks{This research was supported
in part by grants from  ARO (DAAH04-96-1-0193),
AFOSR (F49620-94-1-0415), and
LEQSF (DOD/LEQSF(1996-99)-04).}}}
\author{Xinjia Chen and Kemin Zhou\\
Department of Electrical and Computer Engineering\\
Louisiana State University, Baton Rouge, LA 70803\\
chan@ece.lsu.edu \ \ kemin@ee.lsu.edu}

\date{May 10, 1999}   
\maketitle

\begin{abstract}
In this paper, we consider robust control using randomized algorithms.
We extend the existing order statistics distribution theory
to the general case in which
the distribution of population is not assumed to be continuous and
the order statistics is associated with certain constraints.
In particular, we  derive an inequality on distribution
for related order statistics.
Moreover, we also propose two different approaches in searching
reliable solutions to the robust analysis and
 optimal synthesis problems under constraints.
Furthermore, minimum computational effort is investigated and bounds
 for sample size are derived.
\end{abstract}

\section{Introduction}

It is now well known that many deterministic worst-case
robust analysis and synthesis problems are NP hard, which means that the exact
analysis and synthesis of the corresponding  robust control
problems may be computational demanding \cite{Braatz,Toker}.
On the other hand,  the deterministic worst-case
robustness measures may be quite conservative due to overbounding of
the system uncertainties.
As pointed out  in \cite{KT} by
Khargonekar and  Tikku, the difficulties of
deterministic worst-case robust control problems
are inherent to the problem formulations and a major change of the
paradigm is necessary. An alternative to the deterministic
approach is the probabilistic approach which has been
studied extensively by Stengel and co-workers, see for example,
 \cite{RS,SR}, and references therein.
 Aimed at breaking through the NP-hardness
barrier and reducing the conservativeness of the deterministic
robustness measures, the probabilistic approach  has recently received
a renewed attention in the work by
Barmish and Lagoa \cite{BL},  Barmish, Lagoa, and Tempo \cite{BLT},
Barmish and Polyak \cite{BP}, Khargonekar and  Tikku \cite{KT},
Bai, Tempo, and Fu \cite{bai}, Tempo, Bai,
and Dabbene \cite{TD}, Yoon and Khargonekar \cite{Yoon},
Zhu, Huang and Doyle \cite{XYD}, Chen and Zhou \cite{ChenZ,CZ}
and references therein.

In addition to its low computational complexity,
 the advantages of randomized algorithms can be found in the
flexibility and adaptiveness in dealing with control analysis or synthesis
problems with complicated constraints or in the situation of handling
nonlinearities. The robust control analysis and synthesis problems under
constraints are, in general, very hard to deal with
in the deterministic framework. For example,
it is well-known that a multi-objective control problem involving
mixed $H_2$ and $H_\infty$ objectives are very hard to solve
even though there are elegant solutions to the pure $H_2$ or $H_\infty$
problems \cite{ZDG}.

In this paper,  we first show that most of
the robust control problems can
be formulated as constrained optimal synthesis or
robust analysis problems.
 Since the exact robust analysis or synthesis is,
in general,  impossible, we seek a `reliable' solution
by using randomized algorithms. Roughly speaking,
by  `reliability' we mean how the solution
resulted by randomized algorithms approaches the exact one.
In this paper, we measure the degree of `reliability'
in terms of accuracy $1-\varepsilon$ and confidence level $1-\delta$.
Actually, terminologies like `accuracy' and `confidence level'
have been used in \cite{TD} and \cite{KT} where accuracy $1-\varepsilon$
is referred as an upper bound of the absolute {\it volume} of a subset
of parameter space ${\bf Q}$.  However, in this paper,
we emphasis that the accuracy $1-\varepsilon$ is an upper bound for
the {\it ratio} of {\it volume} of the {\it constrained subset} $
{\bf Q}_{{\bf C}}:=\left\{\;{\rm constraint \ set} \;{\bf C}\; {\rm holds},\;q \in{\bf Q}\right\}$
with respect to the volume of parameter space ${\bf Q}$.
For example,  when estimating the minimum of a quantity $u(q)$ over
$
{\bf Q}_{{\bf C}},
$
the ratio may be
$
\frac{{\rm volume}\;{\rm of}\;\left\{u(q) \geq {\hat{u}}_{min},\;q \in
{\bf Q}_{{\bf C}}\right\}}{{\rm volume} \;{\rm of}\;{\bf Q}_{\bf C}}
$
where ${\hat{u}}_{min}$ is an estimate
resulted by randomized algorithms for quantity $u(q)$.
We can see that the ratio of volume is a better
indicator of the `reliability' than the absolute
${\rm volume}\;{\rm of}\;\left\{u(q) \geq {\hat{u}}_{min},\;q \in {\bf Q}_{{\bf C}}\right\}$.

Based on this measure of `reliability',
we propose two different approaches aimed at
seeking a solution to the robust analysis or
optimal synthesis problem with a certain a priori specified
degree of `reliability'.  One is the {\bf direct approach}.
The key issue is to determine the number of samples
needed to be generated from the parameter space ${\bf Q}$ for a given
reliability measure.
Actually, Khargonekar and Tikku  in
\cite{KT} have applied similar approach to stability margin problem,
though the
measure of `reliability' is in terms of the absolute volume.
In that paper, a sufficient condition is derived on the sample size required
to come up with a `reliable' estimate of the robust
stability margin (See Theorem 3.3 in \cite{KT}).
In this paper, we also derive the bound of sample size and give the sufficient and
necessary condition for the existence of minimum
distribution-free samples size.
 Our result shows that, the bound of sample size necessarily involves
 $\rho:={\rm volume}\;{\rm of}\;{\bf Q}_{{\bf C}}$.
 Thus estimating $\rho$ becomes essential.
 Unfortunately, estimating $\rho$ is time-consuming and
 the resulted sample size is not accurate.
 To overcome this difficulty,  we propose
 and strongly advocate another approach---the {\bf indirect approach}.
 The key issue is to determine the {\it constrained sample size}, which is
 the number of samples needed that fall into
 the constrained subset ${\bf Q}_{\bf C}$.
 We derive bounds of constrained sample size
 and give the sufficient and
necessary condition for the existence of minimum distribution-free
constrained samples size. The bounds do not involve $\rho$
and can be computed exactly. This result makes
it possible to  obtain a reliable solution without estimating
the volume of the constrained parameter subset
${\bf Q_C}$.

This paper is organized as follows. Section 2 presents the
problem formulation and motivations. In Section 3, we derive the exact
distribution of related order statistics
 without the continuity assumption.
Distribution free tolerance interval and
estimation of quantity range
is discussed in Section 4.  Section 5 gives the minimum sample size
under various assumptions.

\section{Preliminary and Problem Formulation}
Let $q={[q_{1}\; \cdots \;q_{n}]}^{T}$ be a vector of
a control system's parameters,
bounded in a compact set ${\bf Q}$, i.e., $q \in {\bf Q}$.
Let ${\bf C}$ be a set of constraints that $q$
must satisfy. Define the {\it constrained subset} of ${\bf Q}$ by
$
{\bf Q}_{{\bf C}}:=\left\{\;{\bf C}\; {\rm holds},\;q \in{\bf Q}\right\}$.
Let $u(q)$ denote a performance index
function. In many applications,
we are concerned with a performance index function
$u(q)$ of the system under the set of constraints ${\bf C}$.
It is natural to ask the following questions:

\begin{itemize}

\item What is $\min_{{\bf Q}_{\bf C}}u(q)$
(or $\max_{{\bf Q}_{\bf C}}u(q)$)?

\item What is the value of $q$ at which $u(q)$ achieves
 $\min_{{\bf Q}_{{\bf C}}}u(q)$ (or $\max_{{\bf Q}_{{\bf C}}}u(q)$)?

\end{itemize}

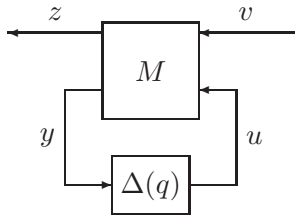
\begin{figure}[htb] 
\centering{
\setlength{\unitlength}{0.050in}%
\begin{picture}(30,21)
\put(10,10){\framebox(10,10){$M$}}       
\put(11,0){\framebox(8,6){$\Delta(q)$}}          
\put(5,21){\makebox(0,0){$z$}}           
\put(5,8){\makebox(0,0)[r]{$y$}}         
\put(25,21){\makebox(0,0){$v$}}          
\put(25,8){\makebox(0,0)[l]{$u$}}        
\put(10,19){\vector(-1,0){10}}
\put(30,19){\vector(-1,0){10}}
\put(24,13){\vector(-1,0){4}}
\put(6,13){\line(1,0){4}}
\put(19,3){\line(1,0){5}}
\put(6,3){\vector(1,0){5}}
\put(6,3){\line(0,1){10}}
\put(24,3){\line(0,1){10}}
\end{picture}
}
\caption{Uncertain System}
\label{fig_a}
\end{figure}

Consider, for example, an uncertain system shown in Figure ~\ref{fig_a}.
Denote the transfer function from $v$ to $z$ by $T_{zv}$ and suppose  that
$T_{zv}$ has the following state space realization
$
T_{zv} = \left[
\begin{array}{c|c}
A (q) & B (q) \\
\hline C(q) & D(q) \end{array} \right].
$
We can now consider several robustness problems:

\begin{itemize}

\item Robust stability: Let ${\bf Q}_{\bf C}= {\bf Q}$ and
$u(q) := \max_{i} {\rm Re} \, \lambda_i  (A(q))$
where $\lambda_i(A)$ denotes the $i$-th eigenvalue of $A$. Then
the system is robustly stable if $\max_{q \in {\bf Q_C}} \  u(q)<0.$

\item Stability margin:
Assume that $\Delta(q)$ belongs to the class of allowable
perturbations ${\bf \Delta}$ which has a certain block structure.
For a given real number $\gamma$, let ${\bf \Delta}_{\gamma}$ denote the subset of perturbations in ${\bf \Delta}$ with size at most $\gamma$, i.e.,
$
{\bf \Delta}_{\gamma}:=\left\{\Delta(q) \in {\bf \Delta}\;:\;\bar{\sigma}(\Delta(q)) \leq \gamma\right\}.
$
The {\it robustness measure} $\gamma_{opt}$  is defined as the smallest allowable perturbation that destabilizes the feedback interconnection.  Let $\gamma_{0}$ be an upper bound for $\gamma_{opt}$.  Define parameter space ${\bf Q}$ by
$
{\bf Q}:=\left\{q\;:\;\Delta(q) \in {\bf \Delta}_{{\gamma}_{0}}\right\}
$
and constrained subset ${\bf Q}_{\bf C}$ by
$
{\bf Q}_{\bf C}:=\{q:q \in {\bf Q} \;{\rm and}\; A(q)
\;{\rm is}\; \mbox{\rm unstable}\}.
$
Let
$
u(q):=\bar{\sigma}(\Delta(q)).
$
It follows that the stability margin problem is equivalent to computing
$
\gamma_{opt}=\min_{{\bf Q}_{\bf C}}u(q).
$
\item Robust performance: Suppose  $A(q)$ is stable for all $q \in {\bf Q}$.
Define
$
u(q):={||T_{zv}||}_{\infty}.
$
Then the robust performance problem is to determine if
$
\max_{q \in {\bf Q}} u(q) \leq \gamma
$
is satisfied for some prespecified $\gamma >0$.

\item Performance range: Let ${\bf Q_C} \subseteq {\bf Q}$ be a given set of
parameters such that $A(q)$ is stable for all $q \in {\bf Q_C}$.
Define again $u(q):={||T_{zv}||}_{\infty}$. Then
the problem of determining  the range of the system's
$H_\infty$  performance level can be formulated as finding
$
\min_{q \in {\bf Q_C}} u(q) \ \ {\rm and} \ \
\max_{q \in {\bf Q_C}} u(q)$.
\end{itemize}

As another example, consider the problem of designing a controller $K(q)$ for
an uncertain system $P(s)$.  Suppose that
$q$ is a vector of controller parameters to be designed
and that the controller is connected with $P(s)$ in a lower LFT setup.
Let the transfer function of the whole system be denoted as
$F_l(P(s),K(q))$.  Suppose that $F_l(P(s),K(q))$ has the following
state space realization
$
F_l(P(s),K(q)) = \left[
\begin{array}{c|c}
A_s (q) & B_s (q) \\
\hline C_s(q) & D_s(q) \end{array} \right]$.
Then we can formulate the problem as a constrained optimal
synthesis problem by defining a performance index
$
u(q):={||F_l(P(s),K(q))||}_{\infty}
$
and restricting parameter $q$ to
$
{\bf Q}_{\bf C} := \left\{\; \max_{i} {\rm Re} \, \lambda_i (A_s(q))< -\alpha, \;q \in {\bf Q}\right\}
$
where $\alpha > 0 $ is not too small for a stability margin.
Then the $H_\infty$ design problem is to
 determine a vector of parameters achieving
$
\min_{q \in {\bf Q_C}}  u(q)$.

\subsection{A Measure of Reliability}

Since the exact solution to the analysis or synthesis problem is impossible.
Measuring how the solution resulted by the randomized algorithm approaches
the exact one becomes essential.
We shall first introduce the concept of {\it volume} proposed  in \cite{KT}.
Let $w(q)$ denote the cumulative distribution function of $q$.  For a subset ${\bf U} \subseteq {\bf Q}$, the
volume of ${\bf U}$, denoted by ${vol}_{w}\{{\bf U}\}$, is defined by
$
{vol}_{w}\{{\bf U}\}:=\int_{q\in {\bf U}} dw(q)$.  Define
$
\rho:=\frac{ {vol}_{w}\left\{{\bf Q}_{\bf C}\right\}}{{vol}_{w}
\left\{{\bf Q}\right\} }$.  Then it follows that $
\rho ={vol}_{w}\left\{{\bf Q}_{\bf C}\right\}
$ since ${vol}_{w}({\bf Q})=1$.  We assume throughout this paper
that $u(q)$ is a measurable function of $q$
and that $vol_{w}\{{\bf Q}_{\bf C}\} > 0$. We also assume
throughout this paper
that $\varepsilon,\delta \in(0,1)$.
Let ${\hat{u}}_{min}$ and  ${\hat{u}}_{max}$ be the estimates
of $\min_{q \in {\bf Q_C}}  u(q)$ and
$\max_{q \in {\bf Q_C}}  u(q)$ respectively.
Note that ${\hat{u}}_{max}$ and ${\hat{u}}_{min}$
are random variables resulted by randomized algorithms.
A reliable estimate of ${\hat{u}}_{min}$ should guarantee
$
{\rm Pr}\left\{\frac{vol_{w}\left\{u(q) \geq {\hat{u}}_{min},\;q \in
{\bf Q}_{{\bf C}}\right\}}{vol_{w}\{{\bf Q}_{\bf C}\}}
\geq 1-\varepsilon\right\} \geq 1-\delta
$
for a small $\varepsilon $ and a small $\delta $.  Similarly,
a reliable estimate of ${\hat{u}}_{max}$ should guarantee
$
{\rm Pr}\left\{
\frac{vol_{w}\left\{u(q) \leq {\hat{u}}_{max},\;q \in {\bf Q}_{\bf C}\right\}}
{vol_{w}\{{\bf Q}_{\bf C}\}} \geq 1-\varepsilon\right\} \geq 1-\delta
$
for a small $\varepsilon $ and a small $\delta$.

\subsection{Two Different Approaches}
\begin{itemize}
\item {\bf Indirect Approach}
Generate i.i.d. samples $q^{i}$ for $q$ by the same distribution
function $w(q)$.  Continue the sampling process until we obtain
$N_{c}$ observations of $q$ which belong to ${\bf Q}_{\bf C}$. Let
$L$ be the number of i.i.d. experiments when this sampling process
is terminated.  Then $L$ is a random number with distribution
satisfying $\sum_{l=N_{c}}^{\infty}{\rm Pr}\left\{L=l\right\}=1$ and
we can show that $E[L]=\frac{N_{c}} {\rho}$.  Let the observations
which belong to ${\bf Q}_{\bf C}$ be denoted as
$q_{c}^{i},\;\;i=1,\cdots,N_{c}$. Define order statistics
${\hat{u}}_{i},\; i=1,\cdots,N_{c}$ as the $i$th smallest one of the
set of observations $\{u(q_{c}^{i}):\; i=1,\cdots,N_{c}\}$, i.e.,
${\hat{u}}_{1} \leq \cdots \leq {\hat{u}}_{N_{c}}$. Obviously, it is
reasonable to take ${\hat{u}}_{1}$ as an estimate for $\min_{{\bf
Q}_{{\bf C}}}u(q)$ and ${\hat{u}}_{N_{c}}$ as an estimate for
$\max_{{\bf Q}_{{\bf C}}}u(q)$ if $N_{c}$ is sufficiently large.
Henceforth, we need to know $N_{c}$ which guarantees $ {\rm
Pr}\left\{\frac{vol_{w}\left\{u(q) \geq {\hat{u}}_{1} ,\;q \in {\bf
Q}_{\bf C}\right\}}{vol_{w}\{{\bf Q}_{\bf C}\}} \geq
1-\varepsilon\right\} \geq 1-\delta $ and (or) $ {\rm Pr}\left\{
\frac{vol_{w}\left\{u(q) \leq {\hat{u}}_{N_{c}},\;q \in {\bf Q}_{\bf
C}\right\}} {vol_{w}\{{\bf Q}_{\bf C}\}} \geq 1-\varepsilon\right\}
\geq 1-\delta$. We call $N_{c}$ {\it constrained sample size}.
\item {\bf Direct Approach}
Let $q^{1},\cdots,q^{N}$ be $N$ i.i.d samples generated by
the same distribution function $w(q)$. Define
${\bf S}:=\left\{q^{1},\cdots,q^{N}\right\}\bigcap{\bf Q}_{{\bf C}}$.
Let $M$ be the number of the elements in ${\bf S}$.
Then $M$ is a random number.
If $M \geq 1$ we denote the
elements of ${\bf S}$ as $q_{c}^{i},\;\;i=1,\cdots,M$.
Define order statistics ${\hat{u}}_{i},\; i=1,\cdots,M$
as the $i$th
smallest one of the set of observations $\{u(q_{c}^{i}):\; i=1,\cdots,M\}$,
i.e., ${\hat{u}}_{1} \leq \cdots \leq {\hat{u}}_{M}$.
 In particular, let ${\hat{u}}_{min}={\hat{u}}_{1}$ and
 ${\hat{u}}_{max}={\hat{u}}_{M}$.  We need to know $N$
which guarantees
$
{\rm Pr}\left\{\frac{vol_{w}\left\{u(q) \geq {\hat{u}}_{min},\;q \in
{\bf Q}_{{\bf C}}\right\}}{vol_{w}\{{\bf Q}_{\bf C}\}}
\geq 1-\varepsilon\right\} \geq 1-\delta
$
and (or)
$
{\rm Pr}\left\{
\frac{vol_{w}\left\{u(q) \leq {\hat{u}}_{max},\;q \in {\bf Q}_{\bf C}\right\}}
{vol_{w}\{{\bf Q}_{\bf C}\}} \geq 1-\varepsilon\right\} \geq 1-\delta$.
We call $N$ {\it global sample size}.
\end{itemize}

\section{Exact Distribution}
Define
$
F_{u}(\gamma) := \frac{vol_{w}\left\{u(q) \leq \gamma
,\;q \in {\bf Q}_{\bf C}\right\}}{vol_{w}\{{\bf Q}_{\bf C}\}}$.
To compute the probabilities involved in Section $2$,
it is important to know
the  associated distribution of any $k$ random
variables $F_{u}({\hat{u}}_{{i}_{1}}),
\;\cdots,\;F_{u}({\hat{u}}_{{i}_{k}}),
\;\;1 \leq {i}_{1} < \cdots < {i}_{k} \leq N_{c},
\;\;1 \leq k \leq N_{c}$ where ${\hat{u}}_{i_{s}},\;s=1,\cdots,k$
is order statistics in the context of the indirect approach.
First, we shall established the following lemma.

\begin{lemma} \label{main3} Let U be a random variable with uniform
distribution over $[0,1]$ and ${\hat{U}}_{n},\;\;
n=1,\cdots,N$ be the order statistics of U,
 i.e., ${\hat{U}}_{1} \leq \cdots \leq {\hat{U}}_{N}$.
Let $0 = t_{0} < t_{1} \cdots < t_{k} \leq 1$.
Define
\[
G_{j_{1}, \cdots ,j_{k}}\left(t_{1},\cdots,t_{k}\right)
:= (1-t_{k})^{N-\sum_{l=1}^{k}j_{l}}\prod_{s=1}^{k} {N-\sum_{l=1}^{s-1}j_{l} \choose j_{s}}
{(t_{s}-t_{s-1})}^{j_{s}}
\]
 and
\[
{\bf I}_ { {i}_{1}, \cdots , {i}_{k} }: =
\left\{(j_{1},\cdots,j_{k}):\;\;i_{s} \leq \sum_{l=1}^{s}j_{l} \leq
N,\;\; s=1,\cdots,k\right\}. \]
  Then $ {\rm Pr}\left\{
{\hat{U}}_{i_{1}} \leq t_{1},\cdots, {\hat{U}}_{i_{k}} \leq
t_{k}\right\} =\sum_{(j_{1},\cdots,j_{k}) \in {\bf I}_ { {i}_{1},
\cdots , {i}_{k} }}G_{j_{1}, \cdots
,j_{k}}\left(t_{1},\cdots,t_{k}\right)$.
\end{lemma}
\begin{pf}
Let $j_{s}$ be the number of samples of U which fall
 into $(t_{s-1},t_{s}],\;\;s=1,\cdots,k$. Then the number of
 samples of U which fall into $[0,t_{s}]$ is $\sum_{l=1}^{s}j_{l}$.
 It is easy to see that the event
 $\left\{{\hat{U}}_{i_{s}} \leq t_{s}\right\}$
 is equivalent to event $\left\{i_{s} \leq \sum_{l=1}^{s}j_{l} \leq N\right\}$.
   Furthermore,
  the event
  $
  \left\{ {\hat{U}}_{i_{1}} \leq t_{1},\cdots,{\hat{U}}_{i_{k}} \leq t_{k}\right\}
  $
  is equivalent to the event
  $
  \left\{i_{s} \leq \sum_{l=1}^{s}j_{l} \leq N,\;\;s=1,\cdots,k\right\}$.
  Therefore,
\begin{eqnarray*}
&   & {\rm Pr}\left\{ {\hat{U}}_{i_{1}} \leq t_{1},\cdots,
{\hat{U}}_{i_{k}} \leq t_{k}\right\} = \sum_{(j_{1},\cdots,j_{k}) \in {\bf I}_ { {i}_{1}, \cdots , {i}_{k} }}\;
\prod_{s=1}^{k}{N-\sum_{l=1}^{s-1}j_{l} \choose j_{s}}
{(t_{s}-t_{s-1})}^{j_{s}}\;{(1-t_{k})}^{N-\sum_{l=1}^{k}j_{l}}\\
& = & \sum_{(j_{1},\cdots,j_{k}) \in {\bf I}_ { {i}_{1}, \cdots ,
{i}_{k} }}\;
G_{j_{1}, \cdots ,j_{k}}\left(t_{1},\cdots,t_{k}\right).
\end{eqnarray*}
\end{pf}

\begin{theorem} \label{main1}
Let $0 = t_{0} < t_{1} \leq \cdots \leq t_{k} \leq 1$ and
 $x_{0}=0$, $x_{k+1}=1$, $i_{0}=0$, $i_{k+1}=N+1$.
Define
$
f_{ {i}_{1}, \cdots , {i}_{k} } (x_{1},\cdots,x_{k}):=
\prod_{s=0}^{s=k}\;N!\;\frac
{ {(x_{s+1}-x_{s})}^{i_{s+1}-i_{s}-1} } {(i_{s+1}-i_{s}-1)!}$ and
$
{\bf D}_{ {p}_{1},\cdots , {p}_{k} }
:=\{(x_{1},\cdots,x_{k}):\;\;0 \leq x_{1} \leq
\cdots \leq x_{k},\;x_{s} \leq p_{s},
\;s=1,\cdots,k\}$.  Define
$
F(t_{1},\cdots,t_{k}) := {\rm Pr}\left\{ F_{u}({\hat{u}}_{i_{1}}) < t_{1},
\cdots,
F_{u}({\hat{u}}_{i_{k}}) < t_{k}\right\}
$
and  $\tau_{s}:
=\sup_{ \{x:\;F_{u}(x) <t_{s} \} }F_{u}(x),\;\;s=1,\cdots,k$.
Then
\[
F(t_{1},\cdots,t_{k}) =
\int_{{\bf D}_{ {\tau}_{1},\cdots,{\tau}_{k} } }
f_{ {i}_{1}, \cdots , {i}_{k} } (x_{1},\cdots,x_{k})\;
dx_{1} \cdots dx_{k} \leq \int_{{\bf D}_{ {t}_{1}, \cdots ,
{t}_{k} }}f_{ {i}_{1}, \cdots , {i}_{k} }
(x_{1},\cdots,x_{k})\;dx_{1} \cdots dx_{k}
\]
and the last equality holds if and only if
$
\exists x_s^* \;{\rm such\; that} \;
{\rm Pr} \{u(q)  < x_s^*\;|\; q \in {\bf Q}_{\bf C} \} = t_s
,\;\;s=1,\cdots,k$.
\end{theorem}

\begin{pf}
Define $\alpha_{0}:= -\infty$ and
${\alpha}_{s} := \sup\left\{x:\;F_{u}(x)
< t_{s}\right\} ,\;\;{\alpha}_{s}^{-}:= \alpha_{s}-{\epsilon},
\;\;s=1,\cdots,k$
where $\epsilon > 0$ can be arbitrary small.
Let $\phi_{s}:= F_{u}({\alpha}_{s}^{-}),\;\;s=1,\cdots,k$. We can show that
${\phi}_{l} < {\phi}_{s}$ if
${\alpha}_{l} < {\alpha}_{s},\; 1 \leq l < s \leq k$.
  In fact, if this is not true, we have
${\phi}_{l} = {\phi}_{s}$. Because $\epsilon$ can be arbitrarily small,
we have $\alpha_{s}^{-} \in (\alpha_{l},\;\alpha_{s})$.  Notice that
${\alpha}_{l} = \min\left\{x:\;F_{u}(x) \geq t_{l} \right\}$
, we have
$t_{l} \leq \phi_{s} = \phi_{l}$.  On the other hand,
by definition we know that
$\alpha_{l}^{-} \in \left\{x:\;F_{u}(x)
< t_{l}\right\}$ and thus $\phi_{l}=F_{u}(\alpha_{l}^{-}) < t_{l}$,
which is a contradiction.  Notice that $F_{u}(\gamma)$ is nondecreasing and right-continuous, we have
${\alpha}_{1} \leq \cdots \leq {\alpha}_{k}$
and $0 \leq {\phi}_{1} \leq \cdots \leq {\phi}_{k} \leq 1
$
and that event
$\left\{F_{u}({\hat{u}}_{i_{s}}) < t_{s}|\;L=l\right\}$
is equivalent to the event
$\left\{{\hat{u}}_{i_{s}} < {\alpha}_{s}|\;L=l\right\}$.
Furthermore, event
$
\left\{
F_{u}({\hat{u}}_{i_{1}}) < t_{1},\;\cdots,
\;F_{u}({\hat{u}}_{i_{k}}) < t_{k}|\;L=l\right\}
$
is equivalent to event
$\left\{{\hat{u}}_{i_{1}} < {\alpha}_{1},\;
\cdots,\;
{\hat{u}}_{i_{k}} < {\alpha}_{k}|\;L=l\right\}
$
which is defined by $k$ constraints
${\hat{u}}_{i_{s}} < {\alpha}_{s}, \;s=1,\cdots,k$.
For every $l <k$, delete constraint ${\hat{u}}_{i_{l}} < {\alpha}_{l}$
if there exists $s >l$ such that ${\alpha}_{s}={\alpha}_{l}$.  Let the remaining constraints
be ${\hat{u}}_{i^{'}_{s}} < {\alpha}^{'}_{s}, \;s=1,\cdots,k^{'}$ where
${\alpha}^{'}_{1} < \cdots < {\alpha}^{'}_{k^{'}}$.
Since all
constraints deleted are actually redundant, it follows that
event $\left\{{\hat{u}}_{i_{1}} < {\alpha}_{1},\;
\cdots,\;{\hat{u}}_{i_{k}} < {\alpha}_{k}|\;L=l\right\}
$ is equivalent to event
$
\left\{{\hat{u}}_{i^{'}_{1}} < {\alpha}^{'}_{1},\;
\cdots,\;
{\hat{u}}_{i^{'}_{k^{'}}} < {\alpha}^{'}_{k^{'}}|\;L=l\right\}$.  Now let $j_{s}$ be the number of observations $u(q_{c}^{i})$
which fall into $[{\alpha}^{'}_{s-1},{\alpha}^{'}_{s}),\;\;s=1,\cdots,k^{'}$.
 Then the number of observations $u(q_{c}^{i})$ which fall into
 $(-\infty,{\alpha}^{'}_{s})$ is $\sum_{l=1}^{s}j_{l}$.
 It is easy to see that the event $\left\{{\hat{u}}_{i^{'}_{s}} <
 {\alpha}^{'}_{s}|\;L=l\right\}$ is equivalent to the event
 $\left\{i^{'}_{s} \leq \sum_{l=1}^{s}j_{l} \leq N\right\}$.
 Furthermore, the event $\left\{{\hat{u}}_{i^{'}_{1}} < {\alpha}^{'}_{1},\;
 \cdots,\;
 {\hat{u}}_{i^{'}_{k^{'}}} < {\alpha}^{'}_{k^{'}}|\;L=l\right\}$
 is equivalent to event
$
\left\{i^{'}_{s} \leq \sum_{l=1}^{s}j_{l} \leq N,
\;\;s=1,\cdots,k^{'}\right\}$.
Therefore
\begin{eqnarray*}&   & {\rm Pr}\left\{ F_{u}({\hat{u}}_{i_{1}}) <
t_{1},\;\cdots,\;
F_{u}({\hat{u}}_{i_{k}}) < t_{k}|\;L=l\right\}\\
& = & {\rm Pr}\left\{{\hat{u}}_{i_{1}} < {\alpha}_{1},\;\cdots,\;
{\hat{u}}_{i_{k}} < {\alpha}_{k}|\;L=l\right\} =
 {\rm Pr}\left\{{\hat{u}}_{i^{'}_{1}} < {\alpha}^{'}_{1},\;
\cdots,\;{\hat{u}}_{i^{'}_{k^{'}}} < {\alpha}^{'}_{k^{'}}|\;L=l\right\}\\
& = & \sum_{(j_{1},\cdots,j_{k^{'}})
\in {\bf I}_ { {i}^{'}_{1}, \cdots , {i}^{'}_{k^{'}} }}
\prod_{s=1}^{k^{'}}{N-\sum_{l=1}^{s-1}j_{l} \choose j_{s}}
{[F_{u}({{\alpha}^{'}_{s}}^{-})-
F_{u}({{\alpha}^{'}_{s-1}}
^{-})]}^{j_{s}}\;{[1-
F_{u}({{\alpha}^{'}_{k^{'}}}^{-})
]}^{N-\sum_{l=1}^{k^{'}}j_{l}}\\
& = & \sum_{(j_{1},\cdots,j_{k^{'}})
\in {\bf I}_ { {i}^{'}_{1}, \cdots , {i}^{'}_{k^{'}} }}
G_{j_{1}, \cdots ,j_{k^{'}}}\left(\phi^{'}_{1},
\cdots,\phi^{'}_{k^{'}}\right).
\end{eqnarray*}
Now consider event
$\left\{ {\hat{U}}_{i_{1}} \leq \phi_{1},\;\cdots,\;{\hat{U}}_{i_{k}} \leq \phi_{k}\right\}$.
For every $l <k$, delete constraint
${\hat{U}}_{i_{l}} \leq {\phi}_{l}$
if there exists $s >l$ such that ${\phi}_{s}={\phi}_{l}$.
Notice that $\phi_{s}=F_{u}({\alpha}_{s}^{-})$ and ${\phi}_{l} < {\phi}_{s}$ if
${\alpha}_{l} < {\alpha}_{s},\; 1 \leq l < s \leq k$, the remaining constraints must
be ${\hat{U}}_{i^{'}_{s}} \leq {\phi}^{'}_{s}, \;s=1,\cdots,k^{'}$ where
${\phi}^{'}_{s} = F_{u}({{\alpha}^{'}_{s}}^{-}), \;s=1,\cdots,k^{'}
$ and $\phi^{'}_{1} < \cdots < \phi^{'}_{k^{'}}$.  Since all
constraints deleted are actually redundant, it follows that event
$\left\{{\hat{U}}_{i_{1}} \leq {\phi}_{1},\;\cdots,\;
{\hat{U}}_{i_{k}} \leq {\phi}_{k}\right\}
$ is equivalent to event
$\left\{{\hat{U}}_{i^{'}_{1}} \leq {\phi}^{'}_{1},\;
\cdots,\;
{\hat{U}}_{i^{'}_{k^{'}}} \leq {\phi}^{'}_{k^{'}}\right\}$.
By Theorem $2.2.3$ in \cite{david} and Lemma~\ref{main3}
 \begin{eqnarray*}
&   & \int_{{\bf D}_{ {\phi}_{1},\cdots,{\phi}_{k} } }
f_{ {i}_{1}, \cdots , {i}_{k} }
(x_{1},\cdots,x_{k})\;dx_{1} \cdots dx_{k} =
 {\rm Pr} \left\{{\hat{U}}_{i_{1}} \leq
{\phi}_{1},\;\cdots,\;{\hat{U}}_{i_{k}} \leq {\phi}_{k}\right\}\\
& = & {\rm Pr}\left\{{\hat{U}}_{i^{'}_{1}} \leq {\phi}^{'}_{1},\;
\cdots,\;
{\hat{U}}_{i^{'}_{k^{'}}} \leq {\phi}^{'}_{k^{'}}\right\} =
 \sum_{(j_{1},\cdots,j_{k^{'}})
\in {\bf I}_ { {i}^{'}_{1}, \cdots , {i}^{'}_{k^{'}} }}
G_{j_{1}, \cdots ,j_{k^{'}}}\left(\phi^{'}_{1},
\cdots,\phi^{'}_{k^{'}}\right).
\end{eqnarray*}
Therefore,
$
{\rm Pr}\left\{ F_{u}({\hat{u}}_{i_{1}}) <
t_{1},\;\cdots,\;
F_{u}({\hat{u}}_{i_{k}}) < t_{k}|\;L=l\right\} =
 \int_{{\bf D}_{ {\phi}_{1},\cdots,{\phi}_{k} } }
f_{ {i}_{1}, \cdots , {i}_{k} }
(x_{1},\cdots,x_{k})\;dx_{1} \cdots dx_{k}$.
It follows that
\begin{eqnarray*}&   & F(t_{1},\cdots,t_{k})
 = {\rm Pr}\left\{ F_{u}({\hat{u}}_{i_{1}}) < t_{1}
 ,\;\cdots
 ,\;F_{u}({\hat{u}}_{i_{k}}) < t_{k}\right\}\\
& = & \sum_{l=N_{c}}^\infty {\rm Pr}\left\{ F_{u}({\hat{u}}_{i_{1}}) < t_{1}
 ,\;\cdots
 ,\;F_{u}({\hat{u}}_{i_{k}}) < t_{k}\;|\;L=l\right\}\;{\rm Pr}\left\{L=l\right\}\\
& = & \sum_{l=N_{c}}^\infty
\int_{{\bf D}_{ {\phi}_{1},\cdots,{\phi}_{k} } }
f_{ {i}_{1}, \cdots , {i}_{k} }
(x_{1},\cdots,x_{k})\;dx_{1} \cdots dx_{k}
\;{\rm Pr}\left\{L=l\right\}
\end{eqnarray*}
Notice that
$
\sum_{l=N_{c}}^\infty {\rm Pr}\left\{L=l\right\} = 1$.
We have
\begin{eqnarray*}&   & F(t_{1},\cdots,t_{k}) =\int_{{\bf D}_{ {\phi}_{1},\cdots,{\phi}_{k} } }
f_{ {i}_{1}, \cdots , {i}_{k} }
(x_{1},\cdots,x_{k})\;dx_{1} \cdots dx_{k} \sum_{l=N_{c}}^\infty
\;{\rm Pr}\left\{L=l\right\}\\
& = & \int_{{\bf D}_{ {\phi}_{1},\cdots,{\phi}_{k} } }
f_{ {i}_{1}, \cdots , {i}_{k} }
(x_{1},\cdots,x_{k})\;dx_{1} \cdots dx_{k}.
\end{eqnarray*}
  By the
definitions of $\tau_{s}$ and $\phi_{s}$, we know that
  ${\bf D}_{ {\tau}_{1},\cdots,{\tau}_{k} }$ is the closure of
  ${\bf D}_{ {\phi}_{1}, {\phi}_{2}, \cdots , {\phi}_{k} }$, i.e.,
   $
  {\bf D}_{ {\tau}_{1},\cdots,{\tau}_{k} }
  = {\bar {\bf D}}_{ {\phi}_{1}, {\phi}_{2}, \cdots , {\phi}_{k} }$
  and that their Lebesgue measures
  are equal.  It follows that
 \[
  F(t_{1},\cdots,t_{k})
  = \int_{ {\bf D}_{ {\tau}_{1},\cdots,{\tau}_{k} } }
  f_{ {i}_{1}, \cdots , {i}_{k} }
  (x_{1},\cdots,x_{k})\;dx_{1} \cdots dx_{k}.
  \]
Notice that $\tau_{s} \leq t_{s},\;\;s=1,\cdots,k$, we have
  $
  {\bf D}_{ {\tau}_{1},\cdots,{\tau}_{k} }
  \subseteq {\bf D}_{ {t}_{1}, \cdots , {t}_{k} }$
  and hence
  \[
  F(t_{1},\cdots,t_{k})
  \leq \int_{{\bf D}_{ {t}_{1}, {t}_{2}, \cdots ,
  {t}_{k} }}f_{ {i}_{1}, \cdots , {i}_{k} }
  (x_{1},\cdots,x_{k})\;dx_{1} \cdots dx_{k},
  \]
  where the equality holds if and only if
  $\tau_{s} = t_{s},\;\;s=1,\cdots,k$,  i.e., $
\exists x_s^* \;{\rm such\; that} \;
{\rm Pr} \{u(q)  < x_s^*\;|\; q \in {\bf Q}_{\bf C} \} = t_s
,\;\;s=1,\cdots,k$.
 \end{pf}

\begin{remark}

For the special case of ${\bf Q}_{\bf C}={\bf Q}$ and
that $F_{u}(.)$ is absolutely continuous, $F(t_{1},\cdots,t_{k})$
can be obtained by combining
{\it Probability Integral Transformation Theorem}
and Theorem $2.2.3$ in \cite{david}.  However, in robust control problem,
the continuity of $F_{u}(.)$ is not
necessarily guaranteed.  For example,
$F_{u}(.)$ is not continuous when uncertain quantity $u(q)$
equals to a constant in an open  set of ${\bf Q}_{\bf C}$.
We can come up with many uncertain systems
in which the continuity assumption for
the distribution of quantity $u(q)$ is not guaranteed.
 Since it is reasonable to assume that
$u(q)$ is measurable, Theorem~\ref{main1} can be applied
in general to tackle these problems without continuity assumption by
a probabilistic approach.  In addition,  Theorem~\ref{main1} can be
applied to investigate the minimum computational effort to come up with a
solution with a certain degree of `reliability' for
robust analysis or optimal synthesis problems under constraints.
\end{remark}

From the proof of Theorem ~\ref{main1}, we can see that
$F(t_{1},\cdots,t_{k})$ is
not related to the knowledge of $L$, thus we have the following corollary.
\begin{corollary} \label{main9}
Let $N_{2} \geq N_{1} \geq N_{c}$.  Then
$
 {\rm Pr}\left\{ F_{u}({\hat{u}}_{i_{1}}) < t_{1},
\;\cdots,\;
F_{u}({\hat{u}}_{i_{k}}) < t_{k}\;\;|\;\;N_{1}
\leq L \leq N_{2}\right\}= F(t_{1},\cdots,t_{k})$.
\end{corollary}

\section{Quantity Range and Distribution-Free Tolerance Intervals}

In robust analysis or synthesis, it is desirable to know function $F_{u}(.)$
because it is actually the distribution function
of quantity $u(q)$ for $q \in {\bf Q}_{\bf C}$.
However, the exact computation of function
$F_{u}(.)$ is in general impossible.
We shall extract as much as possible the
information of $F_{u}(.)$ from observations
$u(q_{c}^{i}),\; i=1,\cdots,N_{c}$.  Let
$
{\cal V}(N_{c},\;i,\;\varepsilon):= \int_{\varepsilon}^{1} \frac{N_{c}!}{(i-1)!(N_{c}-i)!}
x^{i-1}(1-x)^{N_{c}-i}dx
$
for $1 \leq i \leq N_{c}$.

\begin{theorem} \label{cccesu}
$
{\rm Pr}\left\{\frac{vol_{w}\left\{u(q) \geq {\hat{u}}_{m}
,\;q \in {\bf Q}_{\bf C}\right\}}{vol_{w}\{{\bf Q}_{\bf C}\}}
 \geq 1-\varepsilon\right\} \geq 1- {\cal V}(N_{c},\;m,\;\varepsilon)
$
with the equality holds if and only if
$\exists x^*$ such that $F_{u}(x^*) =\varepsilon$.  Moreover,
$
 {\rm Pr}\left\{\frac{vol_{w}\left\{u(q) \leq {\hat{u}}_{m}
,\;q \in {\bf Q}_{\bf C}\right\}}{vol_{w}\{{\bf Q}_{\bf C}\}}
 \geq 1-\varepsilon\right\} \geq 1- {\cal V}(N_{c},\;N_{c}+1-m,\;\varepsilon)
$
with the equality holds if and only if $\exists x^*$
such that ${\rm Pr} \{ u(q)  < x ^*\;|\; q \in {\bf Q}_{\bf C}\} = 1-\varepsilon$.
\end{theorem}

\begin{pf}
Let $v(q)=-u(q)$. Let the cumulative distribution function of $v(q)$ be
$F_{v}(.)$ and define order statistics ${\hat{v}}_{i},\;\;i=1,\cdots,N_{c}$ as the
$i$-th smallest one of the set of observations
$\{v(q_{c}^{i})|\;i=1,\cdots,N_{c}\}$, i.e.,
${\hat{v}}_{1} \leq \cdots \leq {\hat{v}}_{N_{c}}$.
  Obviously, ${\hat{u}}_{m}=-{\hat{v}}_{N_{c}+1-m}$ for any $1 \leq m \leq N_{c}$.
It is also clear that $F_{v}(-x)=1- F_{u}(x^{-})$,
which leads to
$
\sup_{\{ x: F_{v}(x) < 1-\varepsilon\}}F_{v}(x)=1-\varepsilon \Longleftrightarrow
\inf_{\{x: F_{u}(x) > \varepsilon\}}F_{u}(x)=\varepsilon.
$
Apply Theorem ~\ref{main1} to the case of $k=1,\;i_{1}=N_{c}+1-m$, we have
\begin{eqnarray*}
&   & {\rm Pr}\left\{F_{v}({\hat{v}}_{N_{c}+1-m}) < 1-\varepsilon\right\}
=\int_0^{\tau} \frac{N_{c}!}{(N_{c}-m)!(m-1)!}x^{N_{c}-m}
(1-x)^{m-1}dx\\
& \leq & \int_{0}^{1-\varepsilon} \frac{N_{c}!}{(N_{c}-m)!(m-1)!}
x^{N_{c}-m}(1-x)^{m-1}dx={\cal V}(N_{c},\;m,\;\varepsilon)
\end{eqnarray*}
where
$\displaystyle{\tau=\sup_{\{x: F_{v}(x) < 1-\varepsilon\}} F_{v}(x)}$.
Therefore,
\begin{eqnarray*}
&   & {\rm Pr}\left\{\frac{vol_{w}\left\{v(q) \leq {\hat{v}}_{N_{c}+1-m}
,\;q \in {\bf Q}_{\bf C}\right\}}{vol_{w}\{{\bf Q}_{\bf C}\}}
\geq 1-\varepsilon\right\}={\rm Pr}\left\{F_{v}({\hat{v}}_{N_{c}+1-m}) \geq 1-\varepsilon\right\}\\
& = & 1-{\rm Pr}\left\{F_{v}({\hat{v}}_{N_{c}+1-m}) < 1-\varepsilon\right\}
\geq 1-{\cal V}(N_{c},\;m,\;\varepsilon).
\end{eqnarray*}
The equality holds if and only if
$\exists x^*$ such that $F_{u}(x^*) =\varepsilon$ because
such a $x^*$ exists if and only if $\tau=1-\varepsilon$.
It follows that
\begin{eqnarray*}
&   & {\rm Pr}\left\{\frac{vol_{w}\left\{u(q) \geq {\hat{u}}_{m}
,\;q \in {\bf Q}_{\bf C}\right\}}{vol_{w}\{{\bf Q}_{\bf C}\}}
 \geq 1-\varepsilon\right\}\\
& = & {\rm Pr}\left\{\frac{vol_{w}\left\{v(q) \leq {\hat{v}}_{N_{c}+1-m}
,\;q \in {\bf Q}_{\bf C}\right\}}{vol_{w}\{{\bf Q}_{\bf C}\}}
\geq 1-\varepsilon\right\} \geq 1-{\cal V}(N_{c},\;m,\;\varepsilon)
\end{eqnarray*}
with the equality holds if and only if
$\exists x^*$ such that $F_{u}(x^*) =\varepsilon$.

The second part follows by applying
Theorem ~\ref{main1} to the case of $k=1,\;i_{1}=m$.
\end{pf}

\begin{figure}[htb]
\centerline{\psfig{figure=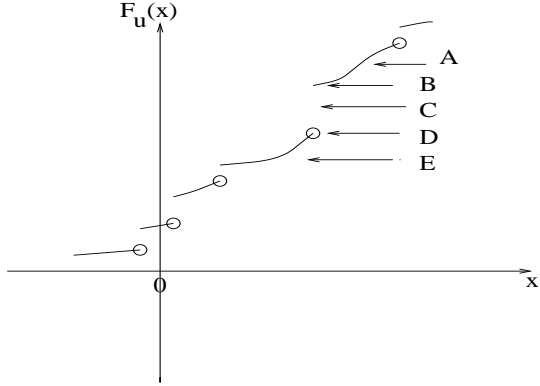,height=2.0in,width=2.8in}}
\caption{
Cases A, B and E guarantee $\exists x^*$ such that $F_u(x^*)=\varepsilon$.
Cases A, D and E guarantee $\exists x^*$ such that
${\rm Pr}\{u(q) < x^*\;|\; q \in {\bf Q}_{\bf C}\}=1-\varepsilon$.
Both conditions are violated in Case C.(The various
magnitude of $\varepsilon$ and $1-\varepsilon$ is
indicated by arrows at different heights.)}
\label{fig}
\end{figure}

It is important to note that the two conditions in Theorem ~\ref{cccesu}
are much weaker than the continuity assumption which requires that
for any $p \in (0,1)$ there exists $x^\star$ such that $F_u(x^\star)=p$.
The difference is visualized in Figure ~\ref{fig}.

In general, it is important to know the probability of a
quantity falling between two arbitrary samples. To that end, we have

\begin{corollary} \label{general}
 Let $1 \leq m < n \leq N_{c}$.
 Suppose $u(q) \neq constant$
 in any open set of ${\bf Q_C}$.  Then
$
{\rm Pr} \{ \frac{vol_{w}\left\{ {\hat{u}}_{m} < u(q) \leq {\hat{u}}_{n},\;q \in {\bf Q}_{\bf C}\right\}}
{vol_{w}\{{\bf Q}_{\bf C}\}} \geq 1-\varepsilon \}=
1-{\cal V}(N_{c},\;N_{c}+1-n+m,\;\varepsilon).
$
\end{corollary}

Since the condition that $u(q) \neq constant$
 in any open set of ${\bf Q_C}$ is equivalent to
 the absolute continuity assumption of $F_{u}(x)$
 (see the proof of Theorem $3.3$ in \cite{KT}),
  the proof of Corollary ~\ref{general}
 can be completed by applying Theorem ~\ref{main1} to the case
 of $k=2$, $i_{1}=m$, $i_{2}=n$ and $F_{u}(x)$ is continuous.

\section{Sample Size}
The important issue of the randomized algorithms to robust analysis or optimal synthesis
is to determine the minimum computational effort required to come up with a solution with a certain degree of `reliability'.  First, we consider
this issue for the indirect approach.

\subsection{Constrained Sample Size}

To estimate $\max_{{\bf Q}_{\bf C}}u(q)$
(or determine parameter $q$ achieving $\max_{{\bf Q}_{\bf C}}u(q)$), we have

\begin{corollary} \label{ccf}
Suppose that
$\exists x^*$ such that ${\rm Pr} \{ u(q)  < x ^*\;|\; q \in {\bf Q}_{\bf C}\} = 1-\varepsilon$. Then
\[
{\rm Pr} \{
\frac{vol_{w}\left\{u(q) \leq {\hat{u}}_{N_{c}},\;q \in {\bf Q}_{\bf C}\right\}}
{vol_{w}\{{\bf Q}_{\bf C}\}} \geq 1-\varepsilon \}
 \geq  1-\delta
\]
 if and only if $N_{c} \geq \frac{\ln{\frac{1}{\delta}}}{\ln{\frac{1}{1-\varepsilon}}}$.
\end{corollary}
It should be noted that  the results in Khargonekar and  Tikku \cite{KT} and Tempo, Bai, and Dabbene \cite{TD} correspond to the sufficient part
of the above Corollary for the special case of ${\bf Q}_{\bf C}={\bf Q}$.

To estimate $\min_{{\bf Q}_{\bf C}}u(q)$
(or determine parameter $q$ achieving $\min_{{\bf Q}_{\bf C}}u(q)$), we have

\begin{corollary} \label{cc}
 Suppose that $\exists x^*$ such that $F_{u}(x^*) =\varepsilon$.  Then
\[
{\rm Pr} \{ \frac{vol_{w}\left\{u(q) \geq {\hat{u}}_{1}
,\;q \in {\bf Q}_{\bf C}\right\}}{vol_{w}\{{\bf Q}_{\bf C}\}}
\geq 1-\varepsilon \} \geq  1-\delta
\]
 if and only if $N_{c} \geq \frac{\ln{\frac{1}{\delta}}}{\ln{\frac{1}{1-\varepsilon}}}$.
\end{corollary}

 To estimate the range of an uncertain quantity with a certain accuracy and confidence level apriori specified, we have the following corollary.

\begin{corollary} \label{ccc}
Suppose that $u(q) \neq constant$
 in any open set of ${\bf Q_C}$.  Then
\[
{\rm Pr}\{ \frac{vol_{w}\left\{ {\hat{u}}_{1} < u(q) \leq {\hat{u}}_{N_{c}}
,\;q \in {\bf Q}_{\bf C}\right\}}
{vol_{w}\{{\bf Q}_{\bf C}\}} \geq 1-\varepsilon \} \geq  1-\delta
\]
if and only if
$\mu(N_{c}) :=
{(1-\varepsilon)}^{N_{c}-1}\left[1+(N_{c}-1)\varepsilon\right]
\leq  \delta$.
\end{corollary}

Now we investigate the computational effort for the direct approach.

\subsection{Global Sample Size }
To estimate $\min_{{\bf Q}_{\bf C}}u(q)$
(or determine parameter $q$ achieving $\min_{{\bf Q}_{\bf C}}u(q)$), we have

\begin{theorem} \label{mm}
Suppose that $\exists x^*$ such that $F_{u}(x^*) =\varepsilon$. Then
\[
{\rm Pr}\left\{\frac{{vol}_{w}\left\{u(q) \geq {\hat{u}}_{min},\;q
\in {\bf Q}_{{\bf C}}\right\}}{{vol}_{w}\left\{{\bf Q}_{\bf
C}\right\}}\geq 1-\varepsilon\right\} \geq 1-\delta \]
 if and only
if $ \displaystyle{ N \geq
\frac{\ln(\frac{1}{\delta})}{\ln(\frac{1}{1-\rho\varepsilon})}}$.
\end{theorem}

\begin{pf}
\begin{eqnarray*}
&   & {\rm Pr}\left\{\frac{{vol}_{w} \left\{u(q) \geq
{\hat{u}}_{min},\;q \in {\bf Q}_{{\bf C}}\right\}}{{vol}_{w}\left\{{\bf
Q}_{\bf C}\right\}}\geq 1-\varepsilon\right\}\\
& = & \sum_{i=0}^{N} {\rm Pr}\left\{M=i\right\}{\rm Pr}\left\{\frac{{vol}_{w}\left\{u(q) \geq {\hat{u}}_{min},\;q \in
{\bf Q}_{{\bf C}}\right\}}{{vol}_{w}\left\{{\bf Q}_{\bf C}\right\}}\geq
1-\varepsilon\;|\;M=i\right\}\\
& = & \sum_{i=0}^{N}{N \choose i}
{\rho}^{i}{(1-\rho)}^{N-i}{\rm Pr}\left\{\frac{{vol}_{w}\left\{u(q) \geq {\hat{u}}_{min},\;q \in
{\bf Q}_{{\bf C}}\right\}}{{vol}_{w}\left\{{\bf Q}_{\bf C}\right\}}\geq
1-\varepsilon\;|\;M=i\right\}.
\end{eqnarray*}
Notice that
$\left\{\frac{{vol}_{w}\left\{u(q) \geq {\hat{u}}_{min},\;q \in
{\bf Q}_{{\bf C}}\right\}}{{vol}_{w}\left\{{\bf Q}_{\bf C}\right\}}\geq
1-\varepsilon\;|\;M=i\right\} \;\Longleftrightarrow\;
\left\{\frac{{vol}_{w}\left\{u(q) \geq {\hat{u}}_{1},\;q \in
{\bf Q}_{{\bf C}}\right\}}{{vol}_{w}\left\{{\bf Q}_{\bf C}\right\}}\geq
1-\varepsilon\;|\;L \leq N\right\}
$
 with $N_{c}=i$ in the context of the indirect
approach.
By Corollary ~\ref{main9}, we know that
\[
{\rm Pr}\left\{\frac{{vol}_{w}\left\{u(q) \geq {\hat{u}}_{1},\;q \in
{\bf Q}_{{\bf C}}\right\}}{{vol}_{w}\left\{{\bf Q}_{\bf C}\right\}}\geq
1-\varepsilon\;|\;L \leq N\right\} = {\rm Pr}\left\{\frac{{vol}_{w}\left\{u(q) \geq {\hat{u}}_{1},\;q \in
{\bf Q}_{{\bf C}}\right\}}{{vol}_{w}\left\{{\bf Q}_{\bf C}\right\}}\geq
1-\varepsilon\right\}.
\]
Apply Theorem ~\ref{cccesu} to the case of $N_{c}=i,\;m=1$, we have
\[
{\rm Pr}\left\{\frac{{vol}_{w}\left\{u(q) \geq {\hat{u}}_{1},\;q \in
{\bf Q}_{{\bf C}}\right\}}{{vol}_{w}\left\{{\bf Q}_{\bf C}\right\}}\geq
1-\varepsilon\right\}
\geq 1-{\cal V}(i,1,\varepsilon)=1-(1-\varepsilon)^{i}
\]
with the equality holds if and only if $\exists x^*$ such that $F_{u}(x^*) =\varepsilon$.  Therefore
\[
{\rm Pr}\left\{\frac{{vol}_{w}\left\{u(q) \geq
{\hat{u}}_{min},\;q \in {\bf Q}_{{\bf C}}\right\}}{{vol}_{w}\left\{{\bf
Q}_{\bf C}\right\}}\geq 1-\varepsilon\right\} \geq
\sum_{i=0}^{N}{N \choose i}
{\rho}^{i}{(1-\rho)}^{N-i}[1-(1-\varepsilon)^{i}]=
1-(1-\varepsilon\rho)^{N}
\]
with the equality holds if and only if
$\exists x^*$ such that $F_{u}(x^*) =\varepsilon$.
Finally, notice that $(1-\varepsilon\rho)^{N} \leq \delta$
if and only if $N \geq \frac{\ln(\frac{1}{\delta})}
{\ln(\frac{1}{1-\rho\varepsilon})}$. This completes the proof.
\end{pf}

It should be noted that sufficiency part of the preceding theorem has been
obtained in \cite{KT} in the context of
estimating robust stability margin.
By the similar argument as that of Theorem ~\ref{mm},
we have the following result for
estimating $\max_{{\bf Q}_{\bf C}}u(q)$
(or determine parameter $q$ achieving $\max_{{\bf Q}_{\bf C}}u(q)$).

\begin{theorem} \label{mmm}
Suppose that
$\exists x^*$ such that ${\rm Pr} \{ u(q)  < x ^*\;|\; q \in {\bf Q}_{\bf C}\} = 1-\varepsilon$.
Then
 \[
{\rm Pr} \{ \frac{{vol}_{w}\left\{u(q) \leq {\hat{u}}_{max},\;q \in
{\bf Q}_{{\bf C}}\right\}}{{vol}_{w}\left\{{\bf Q}_{\bf C}\right\}}\geq
1-\varepsilon \} \geq 1-\delta
\]
 if and only if $\displaystyle{
N \geq \frac{\ln(\frac{1}{\delta})}{\ln(\frac{1}{1-\rho\varepsilon})}.}$
\end{theorem}

To estimate the range of a quantity for the system under a certain constraint ${\bf C}$, we have
\begin{theorem} \label{mmmm}
Suppose $u(q) \neq constant$
in any open set of ${\bf Q_C}$.  Then
\[
{\rm Pr} \{\frac{{vol}_{w}\left\{\;{\hat{u}}_{min} < u(q) \leq {\hat{u}}_{max},\;q \in {\bf Q}_{{\bf C}}\right\}}{{vol}_{w}\left\{{\bf
Q}_{\bf C}\right\}}\geq 1-\varepsilon \} = 1-\mu(N) \geq 1-\delta
\]
if and only if
$\mu(N):={(1-\varepsilon\rho)}^{N-1}[1+(N-1)\varepsilon\rho]\leq \delta$.
\end{theorem}
\begin{pf}

\begin{eqnarray*}
&   & {\rm Pr}\left\{\frac{{vol}_{w}\left\{ {\hat{u}}_{min} < u(q) \leq
{\hat{u}}_{max},\;q \in {\bf Q}_{{\bf C}}\right\}}{{vol}_{w}\left\{{\bf
Q}_{\bf C}\right\}}\geq 1-\varepsilon\right\}\\
& = & \sum_{i=0}^{N} {\rm Pr}\left\{M=i\right\}{\rm Pr}\left\{\frac{{vol}_{w}\left\{{\hat{u}}_{min} < u(q) \leq {\hat{u}}_{max},\;q \in
{\bf Q}_{{\bf C}}\right\}}{{vol}_{w}\left\{{\bf Q}_{\bf C}\right\}}\geq
1-\varepsilon\;|\;M=i\right\}\\
& = & \sum_{i=0}^{N}{N \choose i}
{\rho}^{i}{(1-\rho)}^{N-i}{\rm Pr}\left\{\frac{{vol}_{w}\left\{{\hat{u}}_{min} < u(q) \leq {\hat{u}}_{max},\;q \in
{\bf Q}_{{\bf C}}\right\}}{{vol}_{w}\left\{{\bf Q}_{\bf C}\right\}}\geq
1-\varepsilon\;|\;M=i\right\}.
\end{eqnarray*}
Notice that event $\left\{\frac{{vol}_{w}\left\{ {\hat{u}}_{min} < u(q) \leq {\hat{u}}_{max},\;q \in
{\bf Q}_{{\bf C}}\right\}}{{vol}_{w}\left\{{\bf Q}_{\bf C}\right\}}\geq
1-\varepsilon\;|\;M=i\right\}$ is equivalent to event
\[
\left\{\frac{{vol}_{w}\left\{{\hat{u}}_{1} < u(q) \leq {\hat{u}}_{i},\;q \in
{\bf Q}_{{\bf C}}\right\}}{{vol}_{w}\left\{{\bf Q}_{\bf C}\right\}}\geq
1-\varepsilon\;|\;L \leq N\right\}
\]
 with $N_{c}=i$ in the context of the indirect
approach.

By Corollary ~\ref{main9} and Corollary ~\ref{ccc}, we have
\begin{eqnarray*}
&   & {\rm Pr}\left\{\frac{{vol}_{w}\left\{{\hat{u}}_{1} < u(q) \leq {\hat{u}}_{i},\;q \in
{\bf Q}_{{\bf C}}\right\}}{{vol}_{w}\left\{{\bf Q}_{\bf C}\right\}}\geq
1-\varepsilon\;|\;L \leq N\right\}\\
& = &{\rm Pr}\left\{\frac{{vol}_{w}\left\{{\hat{u}}_{1} < u(q) \leq {\hat{u}}_{i},\;q \in
{\bf Q}_{{\bf C}}\right\}}{{vol}_{w}\left\{{\bf Q}_{\bf C}\right\}}\geq
1-\varepsilon\right\}\\
& = & 1-{(1-\varepsilon)}^{i-1}\left[1+(i-1)\varepsilon\right].
\end{eqnarray*}
Therefore,
\begin{eqnarray*}
&   & {\rm Pr}\left\{\frac{{vol}_{w}\left\{ {\hat{u}}_{min} < u(q) \leq
{\hat{u}}_{max},\;q \in {\bf Q}_{{\bf C}}\right\}}{{vol}_{w}\left\{{\bf
Q}_{\bf C}\right\}}\geq 1-\varepsilon\right\}\\
& = & \sum_{i=0}^{N}{N \choose i}
{\rho}^{i}{(1-\rho)}^{N-i}\left(1-{(1-\varepsilon)}^{i-1}\left[1+(i-1)\varepsilon\right]\right)\\
& = & 1-\sum_{i=0}^{N}{N \choose i}
{\rho}^{i}{(1-\rho)}^{N-i}{(1-\varepsilon)}^{i-1}[1+(i-1)\varepsilon]\\
& = & 1-\frac{1}{1-\varepsilon}\sum_{i=0}^{N}{N \choose i}
{((1-\varepsilon)\rho)}^{i}{(1-\rho)}^{N-i}+\frac{\varepsilon}{1-\varepsilon}\sum_{i=0}^{N}{N \choose i}
{((1-\varepsilon)\rho)}^{i}{(1-\rho)}^{N-i}\\
&   & -N\varepsilon\rho\sum_{i=1}^{N}{N-1 \choose i-1}
{((1-\varepsilon)\rho)}^{i-1}{(1-\rho)}^{N-1-(i-1)}\\
& = & 1-\frac{1}{1-\varepsilon}{(1-\varepsilon\rho)}^{N}
+\frac{\varepsilon}{1-\varepsilon}{(1-\varepsilon\rho)}^{N}
-N\rho\varepsilon{(1-\varepsilon\rho)}^{N-1}\\
& = & 1-{(1-\varepsilon\rho)}^{N-1}[1+(N-1)\varepsilon\rho]\\
& = & 1-\mu(N),
\end{eqnarray*}
which implies that
\[
{\rm Pr}\left\{\frac{{vol}_{w}\left\{{\hat{u}}_{min} < u(q) \leq
{\hat{u}}_{max},\;q \in {\bf Q}_{{\bf C}}\right\}}{{vol}_{w}\left\{{\bf
Q}_{\bf C}\right\}}\geq 1-\varepsilon\right\} \geq 1-\delta
\]
if and only if $\mu(N) \leq \delta$.
\end{pf}

\end{document}